\documentclass{amsart} 
\usepackage{amssymb,latexsym}

\numberwithin{equation}{section}

\newtheorem{theorem}{Theorem}[section]
\newtheorem{proposition}[theorem]{Proposition}
\newtheorem{corollary}[theorem]{Corollary}

\newtheorem{lemma}[theorem]{Lemma}
\newtheorem{example}[theorem]{Example}
\newtheorem{definition}[theorem]{Definition}
\newtheorem{problem}[theorem]{Problem}
\newtheorem{algorithm}[theorem]{Algorithm}

\def\proof{\smallskip\noindent {\bf Proof. }}
\def\endproof{\hfill$\square$\medskip}
\def\endproofmath{\quad\square}

\newcommand{\journal}[1]{{\sl #1}}

\def\CC{\mathbb{C}}

\def\l{\ell}
\def\wnot{w_\mathrm{o}}

\def\cell{X^\circ}

\begin{document}

\title{Recognizing Schubert cells}

\author{Sergey Fomin}
\address{Department of Mathematics, Massachusetts Institute of
  Technology, Cambridge, Massachusetts 02139}
\email{fomin@math.mit.edu}

\author{Andrei Zelevinsky}
\address{\noindent Department of Mathematics, Northeastern University,
  Boston, MA 02115} 
\email{andrei@neu.edu}

\thanks{The authors were supported in part 
by NSF grants \#DMS-9625511 and \#DMS-9700927.
}

\date{\today}




\maketitle

\vspace{-.11in}

\section{Introduction}
\label{sec:intro}

This paper focuses on the properties of Schubert cells 
as quasi-projective subvarieties of a generalized flag variety. 
More specifically, we investigate the problem of distinguishing
between different Schubert cells using 
vanishing patterns of generalized Pl\"ucker coordinates.

\subsection{Formulations of the main problems}
Let $G$ be a simply connected complex semisimple Lie group of rank~$r$
with a fixed Borel subgroup $B$ and a maximal torus $H \subset B$.  
Let $W = {\rm Norm}_G (H)/H$ be the Weyl group of~$G$. 
The generalized flag manifold $G/B$ can be decomposed into the
disjoint union of \emph{Schubert cells}
$\cell_w = (BwB)/B$, for $w\in W$. 

To any weight~$\gamma$ that is $W$-conjugate to some 
fundamental weight of~$G$, one can associate a \emph{generalized
Pl\"ucker coordinate}~$p_\gamma$ on $G/B$ (see \cite{GS} or
Section~\ref{sec:prelim} below). In the case of type~$A_{n-1}$ 
(i.e., $G=SL_n$), 
the $p_\gamma$ are the usual Pl\"ucker coordinates on the flag manifold. 

The closure of a Schubert cell $\cell_w$ is the \emph{Schubert
variety}~$X_w\,$, an irreducible projective subvariety of $G/B$
that can be described as the set of common zeroes of some collection of
generalized Pl\"ucker coordinates~$p_\gamma\,$. 
It is also known 
(see, e.g., Proposition~\ref{prop:1} below) that every Schubert cell 
$\cell_w$ can be defined by specifying vanishing and/or
non-vanishing of some collection of Pl\"ucker coordinates. 

The main two problems studied in this paper are the following.

\begin{problem}
\label{problem:short}
\emph{(Short descriptions of cells)}
Describe a given Schubert cell by as small as possible number of
equations of the form $p_\gamma=0$ and inequalities of the form
$p_\gamma\neq 0$.
\end{problem} 

\begin{problem}
\label{problem:recognition}
\emph{(Cell recognition)}
Suppose a point $x\!\in\! G/B$ is unknown to us, 
but we have access to an 
oracle that answers questions
of the form:~``$p_\gamma (x)\!=\!0$, true or false?'' 
How many such questions are needed to determine the 
Schubert cell $x$ is~in? 
\end{problem}

Problem~\ref{problem:recognition} looks harder than
  Problem~\ref{problem:short}, 
since we do not fix a Schubert cell in advance. 
However, we will 
demonstrate that the complexity of 
the two problems is 
the same: informally speaking, it takes as much time to recognize a cell 
as it takes to describe it. 

Our interest in these problems was originally motivated by their
relevance to the theory of total positivity criteria.
As shown in~\cite{FZ}, 
these criteria take different form in different Bruhat cells~$BwB$, 
so one has to first find out which cell an element $g\in G$ is in. 

\subsection{Overview of the paper}
In Section~\ref{sec:gl3}, we illustrate our problems 
by working out the special case $G=SL_3\,$. 
Section~\ref{sec:prelim} provides the necessary background on
generalized Pl\"ucker coordinates, Bruhat orders, and Schubert
varieties. 

The number of equations of the form
$p_\gamma=0$ needed to define a 
Schubert variety is generally much larger than its codimension. 
In Proposition~\ref{prop:lower-bound}, 
we show that for certain Schubert variety~$X_w$ 
in the flag manifold of type~$A_{n-1}$, 
one needs exponentially many (as a function of~$n$) 
such equations to define it, 
even though ${\rm codim}(X_w)\leq \dim(G/B)=\binom{n}{2}$. 
Given this kind of ``complexity'' of Schubert varieties, 
it may appear surprising that every Schubert \emph{cell} 
actually does have a short description in terms of vanishing or
non-vanishing of certain Pl\"ucker coordinates.   
In Theorem~\ref{th:3}, for the types $A_r$, $B_r$, $C_r$,
and~$G_2$, 
we provide a description of an arbitrary Schubert cell 
$\cell_w$ that only uses ${\rm codim}(X_w)$ equations of the form
$p_\gamma=0$ and at most $r$~inequalities of the form 
$p_\gamma\neq 0$. 
Thus in these cases every Schubert cell is a ``set-theoretic complete
intersection.'' 
Our proof of this property relies on the new concept of an
\emph{economical} linear ordering of fundamental weights. 
For the type~$D$, a description of Schubert cells is slightly more
complicated; see Proposition~\ref{prop:D}. 
This completes our treatment of Problem~\ref{problem:short}. 

In Section~\ref{sec:recognition}, we turn to
Problem~\ref{problem:recognition}.
Our main result is Algorithm~\ref{alg:recognition-general} that
recognizes a Schubert cell $\cell_w$ containing an element~$x$. 
In the cases when an economical ordering exists 
(i.e., for the types $A_r$, $B_r$, $C_r$, and~$G_2$), 
our algorithm ends up examining precisely the same Pl\"ucker
coordinates of~$x$ that appear in Theorem~\ref{th:3}. 
In the case of type~$A_{n-1}\,$, recognizing a
cell requires testing the vanishing of at most $\binom{n}{2}$
Pl\"ucker coordinates. 

In Section~\ref{sec:static}, we 
discuss the 
problem of \emph{cell recognition without feedback},
i.e., the problem of presenting a subset of Pl\"ucker coordinates
whose vanishing pattern determines which cell a point is in.
We show that such a subset must contain
all but a negligible proportion of the Pl\"ucker coordinates.
(Our proof of this result exhibits a surprising connection with
coding theory.)
In Section~\ref{sec:bigrassmannian}, we demonstrate that the situation
changes radically if we only allow \emph{generic} points in each cell. 
With this assumption, knowing the vanishing pattern of
polynomially many Pl\"ucker coordinates (namely, the ones 
corresponding to the \emph{base} of~$W$,
as defined by Lascoux and Sch\"utzenberger~\cite{Las}), 
suffices to recognize a cell.

\subsection{Comments} 
For the purposes of this paper,
all the relevant information about any point on a flag variety 
can be extracted from a finite binary string---the vanishing
pattern of its Pl\"ucker coordinates. 
No explicit description is known for the set of
all possible vanishing patterns. 
For the type~$A$, a combinatorial abstraction of these 
patterns is provided by the notion of a matroid;
for a general Coxeter group, such an abstraction was given by
Gelfand and Serganova~\cite{GS}.
All results of the present paper can be directly extended to 
generalized matroids of~\cite{GS} (irrespective of their
realizability), 
and in fact to a more general combinatorial framework
of ``acceptable'' binary vectors introduced in
Definition~\ref{def:acceptable}. 

\medskip

Note that the ``cell recognition'' problem becomes much simpler if 
its input is an element $gB$ represented by a matrix of~$g$ in some
standard representation of~$G$. 
For instance, if $G=SL_n$, then the Bruhat cell of a given matrix~$g$ 
can be easily determined via Gaussian elimination. 
The reader is referred to~\cite{Gr1}, where an even more general
problem of classifying an 
arbitrary matrix (not necessarily invertible) is solved. 
(This was generalized to the classical series in~\cite{Gr2}.)

\section{Example: $G=SL_3$
}
\label{sec:gl3}

To illustrate our problems,     
let us look at a particular case of type~$A_2$ where $G=SL_3\,$.
In this case, a flag 
$x= (0 \subset F_1\subset F_2\subset F_3=\CC^3) \in G/B$ 
can be represented by a $3 \times 2$ matrix   
whose first column spans $F_1$ and whose first two columns span~$F_2\,$. 
The homogeneous Pl\"ucker coordinates of $x$ are:

\noindent
(1) the matrix entries $p_1$, $p_2$, and $p_3$ in the first column of the
matrix; \nopagebreak

\noindent
(2) the $2\times 2$ minors on the first 2 columns:
$p_{12}$, $p_{13}$, $p_{23}\,$. 

\noindent The complete set of restrictions satisfied by the 6 Pl\"ucker
coordinates consists of: \nopagebreak

\noindent
(a) the \emph{Grassmann-Pl\"ucker relation} 
$p_1 p_{23} - p_2 p_{13} + p_3 p_{12}=0$;\nopagebreak

\noindent
(b) non-degeneracy conditions: $(p_1,p_2,p_3)\neq (0,0,0)$,
 $(p_{12},p_{13},p_{23})\neq (0,0,0)$.

The Weyl group here is the symmetric group~$\mathcal{S}_3\,$,
with generators $s_1 \!=\! (1,2)$ and $s_2 \!=\! (2,3)$
 and relations
$s_1^2=s_2^2=1$ and $\wnot=s_1s_2s_1=s_2s_1s_2$. 
In Table~\ref{table:gl3}, we show which Pl\"ucker coordinates must or
must not vanish on each particular Schubert cell. 
In the table, 
0 means ``vanishes on $\cell_w\,$,''
1 means ``does not vanish anywhere on $\cell_w\,$,''
and the wildcard $*$ means that both zero and nonzero values do
occur.

\begin{table}[ht]

\vspace{.2in}

\begin{tabular}{ll|cccccc|l}
 & &  & & &  & &  & \\[-.15in]
$w$      &  &   $p_1$ &   $p_2$ &   $p_3$ & $p_{12}$ & $p_{13}$ &
$p_{23}$ & 
\qquad $\cell_w$\\[.07in]
\hline
\\[-.25in]
 & &  & & &  & &  & \\
$e$   & 123     &    1  &    0  &    0  &    1   &   0   &   0   &  
$p_3=p_2=p_{13}= 0$  \\[.1in]
$s_1$  & 213    &    $*$  &    1  &    0  &    1   &   0   &   0   &  
$p_{13}= p_{23}   =  0$, $p_2  \neq 0$ \\[.1in]
$s_2$  & 132    &    1  &    0  &    0  &    $*$   &   1   &   0   &   
$p_2 = p_3 = 0$, $p_{13} \neq 0$ \\[.1in]
$s_1 s_2$ & 231 &    $*$  &    1  &    0  &    $*$   &   $*$   &   1   &  
$p_3   =  0$, $p_{23} \neq 0$  \\[.1in]
$s_2 s_1$ & 312 &    $*$  &    $*$  &    1  &    $*$   &   1   &   0   &   
$p_{23}=0$, $ p_3 \neq 0$ \\[.1in]
$\wnot$ & 321   &    $*$  &    $*$  &    1  &    $*$   &   $*$   &   1   &
 $p_3 \neq 0$, $p_{23} \neq 0$ \\[-.13in]
 & &  & & &  & &  & 
\end{tabular}

\vspace{.2in}

\caption{Schubert cells and Pl\"ucker coordinates in type~$A_2$}
\label{table:gl3}
\end{table}

Concerning Problem~\ref{problem:short},
we see that each Schubert cell can be described
in terms of the 4 Pl\"ucker coordinates
$p_2,p_3,p_{13},p_{23}$
(these are exactly the ``bigrassmannian'' coordinates discussed in
Section~\ref{sec:bigrassmannian}).  
Moreover, 3 equations/inequalities suffice to
describe every single cell,
as shown in the last column of Table~\ref{table:gl3}. 

\begin{figure}[ht]
\setlength{\unitlength}{2pt} 
\begin{center}
\begin{picture}(80,105)(0,-5)

\put(0,20){\circle*{1.5}}
\put(0,40){\circle*{1.5}}
\put(0,60){\circle*{1.5}}
\put(30,70){\circle*{1.5}}
\put(40,0){\circle*{1.5}}
\put(40,50){\circle*{1.5}}
\put(40,90){\circle*{1.5}}
\put(50,70){\circle*{1.5}}
\put(80,20){\circle*{1.5}}
\put(80,40){\circle*{1.5}}
\put(80,60){\circle*{1.5}}

\put(0,20){\line(0,1){20}}
\put(0,40){\line(0,1){20}}
\put(80,20){\line(0,1){20}}
\put(80,40){\line(0,1){20}}
\put(40,0){\line(-2,1){40}}
\put(40,0){\line(2,1){40}}
\put(0,20){\line(2,1){80}}
\put(80,20){\line(-2,1){80}}
\put(40,0){\line(0,1){50}}
\put(0,40){\line(1,1){30}}
\put(80,40){\line(-1,1){30}}
\put(0,60){\line(4,3){40}}
\put(80,60){\line(-4,3){40}}
\put(30,70){\line(1,2){10}}
\put(40,50){\line(1,2){10}}
\put(40,50){\line(-1,2){10}}
\put(50,70){\line(-1,2){10}}

\put(40,-6){\makebox(0,0){$0000$}}

\put(-9,20){\makebox(0,0){$1000$}}
\put(-9,40){\makebox(0,0){$1001$}}
\put(-9,60){\makebox(0,0){$1011$}}
\put(89,20){\makebox(0,0){$0010$}}
\put(89,40){\makebox(0,0){$0110$}}
\put(89,60){\makebox(0,0){$1110$}}

\put(47,51){\makebox(0,0){$0101$}}
\put(22,70){\makebox(0,0){$1101$}}
\put(58,70){\makebox(0,0){$0111$}}

\put(40,95){\makebox(0,0){$1111$}}

\put(-3,37){\dashbox{1}(6,26)[t]{}}
\put(77,37){\dashbox{1}(6,26)[t]{}}

\put(-3,17){\dashbox{1}(6,6)[t]{}}
\put(77,17){\dashbox{1}(6,6)[t]{}}
\put(37,-3){\dashbox{1}(6,6)[t]{}}

\put(28,48){\dashbox{1}(24,44)[t]{}}

\put(1,13){\makebox(0,0){$s_1$}}
\put(81,13){\makebox(0,0){$s_2$}}
\put(43,6){\makebox(0,0){$e$}}
\put(0,67){\makebox(0,0){$s_1s_2$}}
\put(80,67){\makebox(0,0){$s_2s_1$}}
\put(57,90){\makebox(0,0){$\wnot$}}

\end{picture}
\end{center}

\vspace{.1in}

\caption{Vanishing patterns of Pl\"ucker coordinates 
$p_2$, $p_3$, $p_{13}$, $p_{23}$
} 
\label{fig:gl3-patterns}
\end{figure}
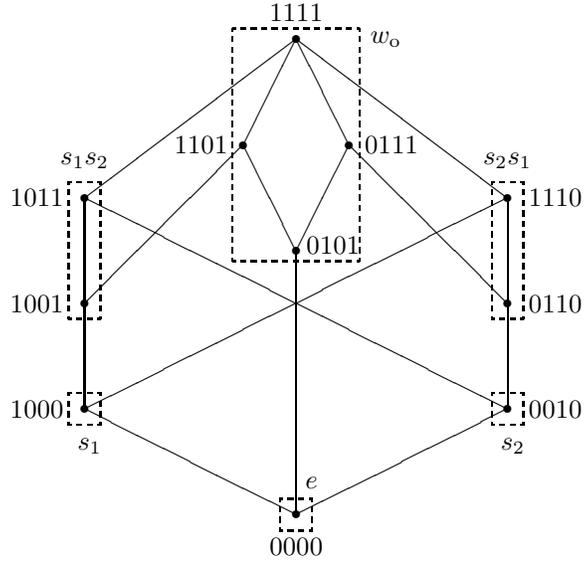

Altogether, there are 11 possible vanishing patterns for
the Pl\"ucker coordinates $p_2,p_3,p_{13},p_{23}\,$. 
The classification of points on the flag variety according to the 
vanishing patterns of these coordinates
provides a refinement of the Schubert cell decomposition. 
In Figure~\ref{fig:gl3-patterns}, we represent this stratification
by a graph (actually, the Hasse diagram of a poset)
whose 11 vertices are labelled by the vanishing patterns 
and whose edges show how the subcells degenerate 
into each other when a condition of the form
$p_I \neq 0$ is replaced by $p_I=0$. 
The dashed boxes 
enclose the subsets making up individual Schubert cells. 
See Section~\ref{sec:bigrassmannian} for further discussion of this
poset.

The Schubert varieties $X_w$ are defined by the 
\emph{equalities} appearing in the last column of Table~\ref{table:gl3}.  
Thus in this case the minimal number of equations
of the form $p_\gamma (x)=0$ that define a Schubert variety $X_w$ as a
subset of $G/B$ is equal to its codimension. 
In general, however, such a  statement is grossly false (see
Section~\ref{sec:lower-bounds-for-the-number-of-equations}).  

\medskip

\begin{figure}[ht]
\setlength{\unitlength}{2.5pt} 
\begin{center}
\begin{picture}(130,65)(-2,-2)

\put(0,0){\circle*{1.5}}
\put(20,0){\circle*{1.5}}
\put(40,0){\circle*{1.5}}
\put(60,0){\circle*{1.5}}
\put(10,20){\circle*{1.5}}
\put(50,20){\circle*{1.5}}
\put(90,20){\circle*{1.5}}
\put(130,20){\circle*{1.5}}
\put(30,40){\circle*{1.5}}
\put(110,40){\circle*{1.5}}
\put(70,60){\circle*{1.5}}

\put(0,0){\line(1,2){10}}
\put(40,0){\line(1,2){10}}
\put(20,0){\line(-1,2){10}}
\put(60,0){\line(-1,2){10}}
\put(10,20){\line(1,1){20}}
\put(90,20){\line(1,1){20}}
\put(50,20){\line(-1,1){20}}
\put(130,20){\line(-1,1){20}}
\put(30,40){\line(2,1){40}}
\put(110,40){\line(-2,1){40}}

\put(50,55){\makebox(0,0){$p_3\!=\!0$}}
\put(90,55){\makebox(0,0){$p_3\!\neq\! 0$}}

\put(16,33){\makebox(0,0){$p_2\!=\!0$}}
\put(94,33){\makebox(0,0){$p_{23}\!=\!0$}}
\put(46,33){\makebox(0,0){$p_2\!\neq\! 0$}}
\put(126,33){\makebox(0,0){$p_{23}\!\neq\! 0$}}

\put(-1,13){\makebox(0,0){$p_{13}\!=\!0$}}
\put(22,13){\makebox(0,0){$p_{13}\!\neq\! 0$}}
\put(39,13){\makebox(0,0){$p_{23}\!=\!0$}}
\put(62,13){\makebox(0,0){$p_{23}\!\neq\! 0$}}

\put(0,-4){\makebox(0,0){123
}}
\put(40,-4){\makebox(0,0){213
}}
\put(20,-4){\makebox(0,0){132
}}
\put(60,-4){\makebox(0,0){231
}}
\put(90,16){\makebox(0,0){312
}}
\put(130,16){\makebox(0,0){321
}}

\end{picture}
\end{center}

\vspace{.1in}

\caption{Cell recognition algorithm for the type $A_2$ 
}
\label{fig:3steps}
\end{figure}
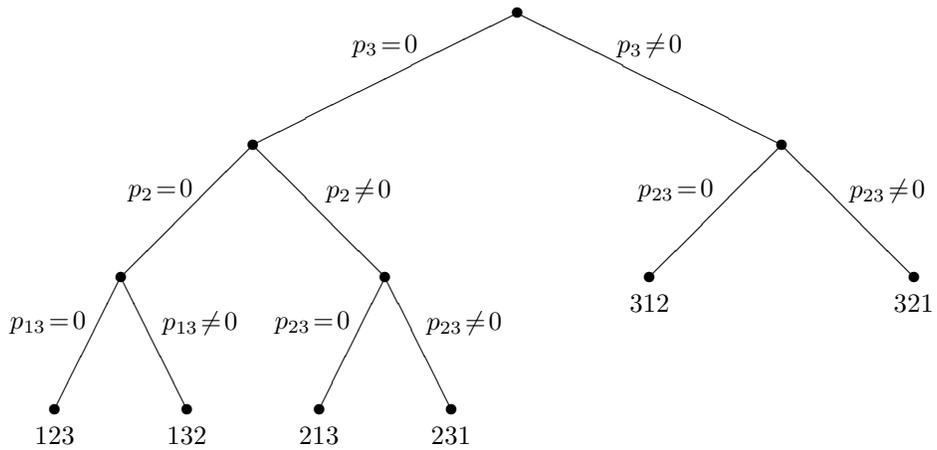

Turning to Problem~\ref{problem:recognition},
the best recognition algorithm is given in Figure~\ref{fig:3steps};
it requires 3 questions. 
Notice that each branch of the tree provides a short description of 
the corresponding Schubert cell.

\section{Preliminaries}
\label{sec:prelim}

In this section, we review basic facts about generalized Pl\"ucker
coordinates, the Bruhat orders, and Schubert varieties.
For general background on these topics, see, e.g.,
\cite[Section~4]{GS}, \cite{bourbaki}, and 
\cite[\S23.3--23.4]{fulton-harris}. 

\subsection{Generalized Pl\"ucker coordinates} 
\label{sec:plucker}

Our approach to this classical subject is similar to the one of Gelfand
and Serganova~\cite[Section~4.2]{GS}. 
Let us fix some linear ordering $\omega_1, \ldots, \omega_r$ of
fundamental weights; the choice of this ordering will later become 
important.  
We will call the weights $\gamma \in W \omega_i$ 
\emph{Pl\"ucker weights of level}~$i$.
Recall that the orbits of fundamental weights are pairwise disjoint, 
so the notion of level is well defined.

Let $V_{\omega_i}$ be the fundamental representation of $G$ with 
highest weight $\omega_i\,$.
For any Pl\"ucker weight~$\gamma$ of level~$i$,
the weight subspace $V_{\omega_i} (\gamma)$ is known to be one-dimensional.
Let us fix an arbitrary nonzero vector 
$v_\gamma \in V_{\omega_i}(\gamma)$ for each such~$\gamma$. 
In particular, $v_{\omega_i}$ is a highest weight vector in
$V_{\omega_i}\,$, and thus an eigenvector for the action of any $b \in
B$; we will write $b v_{\omega_i} = b^{\omega_i} v_{\omega_i}$. 

\begin{definition}
{\rm 
The \emph{generalized Pl\"ucker coordinate}~$p_\gamma$
associated to a Pl\"ucker weight $\gamma$ of level~$i$ 
is defined as follows.
For $g\in G$, let $p_\gamma (g)$ be the coefficient of $v_\gamma$
in the expansion of $g v_{\omega_i}$ into any basis of $V_{\omega_i}$ 
consisting of weight vectors. 
It follows that $p_\gamma (gb) = b^{\omega_i} p_\gamma (g)$ 
for any $g \in G$ and $b \in B$.
Thus we can think of $p_\gamma$ as a global section of the line bundle
on $G/B$ corresponding to the character $b \mapsto b^{\omega_i}$ of~$B$.
It then makes sense to talk about vanishing or non-vanishing of
$p_\gamma$ at any point $x = gB$ of the generalized flag
manifold~$G/B$. 
}\end{definition} 

Although the definition of $p_\gamma$ depends on the choice of
normalization for the vectors $v_\gamma$, this dependence is
not very essential: 
a different choice of normalizations only changes each $p_\gamma$ 
by a nonzero scalar multiple. 
In particular, the set of zeroes of each $p_\gamma$ is a uniquely and
unambiguously defined hypersurface in~$G/B$. 

We note that one natural choice of normalization is the following: 
define $p_\gamma$ as the ``generalized minor''
$\Delta_{\gamma, \omega_i}$, 
in the notation of~\cite[Section~1.4]{FZ}.

For the type $A_{n-1}$, the notion of a Pl\"ucker coordinate specializes 
to the ordinary one (see, e.g.,~\cite{fulton-YT}), as follows. 
Let us use the standard numeration of the fundamental weights, 
so that $V_{\omega_1} = V = \CC^n$ is the defining representation of 
$G = SL_n$, and $V_{\omega_i} = \Lambda^i V$. 
Pl\"ucker weights of level~$i$ are naturally identified with subsets
$I \subset [1,n]$ of cardinality $i$: under this identification,
the weight subspace $V_{\omega_i} (\gamma)$ is the one-dimensional subspace
$\Lambda^i \CC^I \subset \Lambda^i V$. 
The variety $G/B$ is identified with the manifold of all complete
flags $x=(0\subset F_1\subset \cdots \subset F_n = V)$ in~$V$:  
for $x=gB\in G/B$, the subspace $F_i$ is generated by the first $i$
columns of the matrix~$g$.
The Pl\"ucker coordinate $p_I(x)$ is simply the minor of $G$ with the
row set~$I$ and the column set~$[1,i] = \{1, \dots, i\}$.
It follows that  $p_I$ does \emph{not} vanish
at a flag $x$ if and only if $F_i \cap \CC^{[1,n] - I} = \{0\}$.

\subsection{Bruhat orders}
\label{sec:Bruhat}
The Bruhat order can be defined for an arbitrary Coxeter group~$W$.
(Even though it seems to be well established that
the Bruhat order is actually due to Chevalley,
we stick with the traditional terminology to avoid misconceptions.)
Let $S=\{s_1,\dots,s_r\}$ be the set of simple reflections in $W$, and
$\l (w)$ be the length function.
The Bruhat order on $W$ is the transitive closure of the 
following relation:
$w < wt$ for any reflection $t$ (that is, a $W$-conjugate of a simple
reflection) such that $\l (w) <\l (wt)$.

For every subset $J$ of $[1,r]$, 
let $W_J$ denote the parabolic subgroup of $W$
generated by the simple reflections $s_j$ with $j \in J$.
Each coset in $W/W_J$ has a unique representative 
which is minimal with respect to the Bruhat order.  
These representatives are partially ordered by the 
Bruhat order, inducing a partial order on $W/W_J$.
This partial order is also called the Bruhat order on $W/W_J$.

We will be especially interested in the coset spaces modulo 
\emph{maximal} parabolic subgroups $W_{\widehat i} = W_{[1,r] - \{i\}}$.
The following basic result is due to
Deodhar~\cite[Lemma~3.6]{deodhar}.

\begin{lemma}
\label{lem:deodhar}
For $u,v\in W$, we have: $u\leq v$ if and only if 
$u W_{\widehat i} \leq v W_{\widehat i}$ for all $i$.  
\end{lemma}

 From now on we assume that $W$ is the Weyl group associated 
to a semisimple complex Lie group $G$. 
Then the stabilizer of a fundamental weight $\omega_i$
is the maximal parabolic subgroup $W_{\widehat i}$.
Thus the correspondence $w \mapsto w \omega_i$ establishes 
a bijection between the coset space $W/W_{\widehat i}$ and the set
$W \omega_i$ of Pl\"ucker weights of level~$i$. 
This bijection transfers the Bruhat order from
$W/W_{\widehat i}$ to~$W \omega_i$.
Note that if $\gamma$ and $\delta$ are two Pl\"ucker weights 
of the same level, with $\gamma \leq \delta$ with respect to the
Bruhat order, then the weight $\gamma - \delta$ can be expressed as a
sum of simple roots.
The converse statement is true for type $A$ but false in general.  
A counterexample for the type~$B_3$ is given in
\cite[pp.~176--177]{SVZ}; 
see also Deodhar~\cite{deodhar2}
(we thank John Stembridge for providing
this reference). 

The Bruhat order on the Weyl group $W$ also has the following well-known
geometric interpretation  in terms of Schubert cells and Schubert varieties:
\begin{equation}
\label{eq:Schubert-Bruhat}
u \leq v \Longleftrightarrow \cell_u \subset X_v \Longleftrightarrow 
X_u \subset X_v \ .
\end{equation}
A similar interpretation exists for the Bruhat order on any coset space 
$W/W_J$: if $P_J$ is the parabolic subgroup in $G$ corresponding 
to $W_J$ then the correspondence $w \mapsto w P_J$ establishes 
a bijection between $W/W_J$ and $G/P_J$, and we have $u W_J \leq v W_J$
if and only if the ``cell'' $(B u P_J)/P_J$ is contained in the closure 
of $(B v P_J)/P_J$.

To illustrate the above concepts, 
consider the case of type~$A_{n-1}$ where $G=SL_n\,$,
and $W$ is the symmetric group~$\mathcal{S}_n\,$.
We have already seen that Pl\"ucker weights of level $i$
are in natural bijection with the $i$-subsets of~$[1,n]$.
The Bruhat order on the $i$-subsets of $[1,n]$ can be explicitly
described as follows: for two subsets 
$J=\{j_1<\cdots<j_i\}$ and $K=\{k_1<\cdots<k_i\}$, 
we have $J\leq K$ if and only if $j_1\leq k_1,\dots,j_i\leq k_i\,$. 
Lemma~\ref{lem:deodhar} tells that $u\leq v$ in the Bruhat order
if and only if $u([1,i])\leq v([1,i])$ for any~$i$,
in the sense just defined. 
(This is the original Ehresmann's criterion~\cite{ehresmann}.)

\subsection{Set-theoretic description of Schubert varieties}

\nopagebreak[4]

\begin{proposition}
\label{prop:Xw=zeroes}
A point $x \in G/B$ belongs to the Schubert variety $X_w$ if and only
if $p_\gamma (x) = 0$ for any Pl\"ucker weight $\gamma$ 
(say, of level~$i$) such that $\gamma\not\leq w\omega_i$ in the Bruhat
order. 
\end{proposition}

This proposition is well known to experts 
but we were unable to find it explicitly stated in the literature.
It can be deduced from much stronger results in~\cite{Lak,Ram} that
provide a scheme-theoretic description of $X_w$. 
Following the suggestion of Peter Littelmann, we provide a self-contained
proof of Proposition~\ref{prop:Xw=zeroes} which is much more elementary
than the arguments in ~\cite{Lak,Ram}.
We deduce Proposition~\ref{prop:Xw=zeroes} from the following lemma. 

\begin{lemma}
\label{lem:who-vanishes}
A Pl\"ucker coordinate $p_\gamma$ of level~$i$ 
does not identically vanish on the Schubert variety~$X_w$
if and only if $\gamma \leq w\omega_i\,$. 
Also, $p_{u \omega_i}$ vanishes nowhere on $\cell_u\,$. 
\end{lemma}

\proof 
Let us recall the definition of $p_\gamma (gB)$:
up to a nonzero scalar, this is the coefficient of $v_\gamma$
in the expansion of $g v_{\omega_i}$ in any basis of $V_{\omega_i}$ 
consisting of weight vectors. 
Here $V_{\omega_i}$ is the fundamental representation of $G$ with 
highest weight $\omega_i$, and $v_\gamma \in V_{\omega_i}$ is a 
(unique up to a scalar) vector of weight $\gamma$. 
It follows that if $g \in B u B$ and $\gamma=u\omega_i$, 
then $p_\gamma (gB)$ is a nonzero scalar 
multiple of the coefficient of $v_\gamma$ in the expansion of $b
v_\gamma$ for some $b \in B$. 
This coefficient is clearly nonzero which proves the last statement of
the lemma.

Now let $\gamma$ be an arbitrary Pl\"ucker weight of level $i$.
First let us assume that $\gamma \leq w\omega_i$. 
By definition of the Bruhat order, $\gamma = u \omega_i$ for some 
$u \leq w$. 
We have just proved that $p_\gamma$ vanishes nowhere on $\cell_u$. 
But $\cell_u \subset X_w$ by (\ref{eq:Schubert-Bruhat});
therefore $p_\gamma$ does not identically vanish on $X_w$.

It remains to prove the converse statement: if 
$p_\gamma$ does not identically vanish on $X_w$ then $\gamma \leq w\omega_i$.
(The following argument was shown to us by Peter Littelmann;
it closely follows the proof of Proposition~1 in 
Gelfand and Serganova~\cite[Section~5]{GS}.)
Let $P(V_{\omega_i})$ denote the projectivization of the vector space 
$V_{\omega_i}$, and let $[v] \in P(V_{\omega_i})$ denote the projectivization
of a nonzero vector $v \in V_{\omega_i}$. 
Then the stabilizer of $[v_{\omega_i}]$ in $G$ is the maximal parabolic 
subgroup $P_{\widehat i}$, so the map $g \mapsto g [v_{\omega_i}]$
identifies the coset space $G/P_{\widehat i}$ with the orbit
$G [v_{\omega_i}] \subset P(V_{\omega_i})$. 

We shall use the following well known fact: the convex hull of all weights
of the representation $V_{\omega_i}$ is a convex polytope whose vertices
are precisely Pl\"ucker weights of level $i$. 
It follows that, for every Pl\"ucker weight $\gamma$ of level $i$, 
there exists a one-parameter subgroup $\chi: \CC_{\neq 0} \to H$
such that 
\begin{equation}
\label{eq:gamma is extremal}
\lim_{t \to \infty} \chi (t)^\delta / \chi (t)^\gamma = 0
\end{equation}
for any weight $\delta \neq \gamma$ of $V_{\omega_i}$. 

Now everything is ready for concluding the proof.
Suppose $p_\gamma$ does not identically vanish on $X_w$.
By definition, this means that $v_\gamma$ appears with nonzero coefficient
in the expansion of $g v_{\omega_i}$ for some $g \in BwB$.
Using (\ref{eq:gamma is extremal}) we see that 
$$[v_\gamma] = \lim_{t \to \infty} \chi (t) g [v_{\omega_i}] \ .$$
It follows that $[v_\gamma]$ lies in the closure of $(B w B) [v_{\omega_i}]$.
We have $\gamma = u \omega_i$ for some $u \in W$. 
Identifying as above the orbit $G [v_{\omega_i}]$ with the coset space
$G/ P_{\widehat i}$ we conclude that the coset $u P_{\widehat i}$
is contained in the closure of $(B w P_{\widehat i})/ P_{\widehat i}$.
As explained in Section~\ref{sec:Bruhat}, this implies that 
$\gamma = u \omega_i \leq w \omega_i$, and we are done. 
\endproof

We can now complete the proof of Proposition~\ref{prop:Xw=zeroes}.

\proof
Let $X \subset G/B$ denote the variety defined by the equations
$p_\gamma (x) = 0$ for all Pl\"ucker weights $\gamma$ of any level 
$i$ such that $\gamma\not\leq w\omega_i$. 
The inclusion $X_w\subset X$ follows from Lemma~\ref{lem:who-vanishes}. 
Now assume that $x\notin X_w\,$; 
say, $x\in\cell_v$ with $v\not\leq w$. 
By Lemma~\ref{lem:deodhar}, there exists $i$ such that 
$v\omega_i\not\leq w\omega_i\,$. 
By Lemma~\ref{lem:who-vanishes}, $p_{v\omega_i}(x) \neq 0 \,$.
Therefore $x\notin X$, as desired. 
\endproof

\section{Short descriptions of Schubert cells}
\label{sec:short}

This section is devoted to set-theoretic descriptions of Schubert
cells. 
It is well known that $\cell_w$ is defined inside $X_w$ by $r$
inequalities $p_{w \omega_i} \neq 0$.
Combining this with the set-theoretic description of $X_w$ in 
Proposition~\ref{prop:Xw=zeroes}, we obtain a set-theoretic description 
of~$\cell_w$. 
However, the following proposition shows that we can do better.

\begin{proposition}
\label{prop:1}
An element $x\in G/B$ belongs to a Schubert cell $\cell_w$ if and
only~if, for every $i\in[1,r]$, the following conditions hold: 
\smallskip

{\rm(1)}\quad $p_{w \omega_i} (x) \neq 0$;
\smallskip

{\rm(2)}\quad  $p_{\gamma} (x) = 0$ for all 
$\gamma \in  w W_{[i,r]} \omega_i$ such that $\gamma > w \omega_i\,$.
\end{proposition}

\proof
In view of Lemma~\ref{lem:who-vanishes},
these conditions are certainly neccesary.
Let us prove that they are also sufficient.
Suppose that (1)--(2) hold, and let $x\in\cell_u\,$.
Our goal is to show that $u=w$.
First, by (1) and Lemma~\ref{lem:who-vanishes}, we have 
$w \omega_i \leq u \omega_i$ for all $i$, hence $w \leq u$ 
by Lemma~\ref{lem:deodhar}. 
Now suppose that $w < u$. 
Then at least one of the inequalities $w \omega_i \leq u \omega_i$ 
is strict; take the minimal index $i$ such that $w \omega_i < u
\omega_i$.  
The equalities $w \omega_j = u \omega_j$ for $j < i$ imply that
$w^{-1} u \in \bigcap_{j=1}^{i-1} W_{\widehat j}$.
Using 
the equality 
$W_{J_1} \cap \cdots \cap W_{J_k} = W_{J_1 \cap \cdots \cap J_k}$
valid in any Coxeter group (see \cite{bourbaki}),
we conclude that
\begin{equation}
\label{eq:intersection-of-parabolics}
\bigcap_{j=1}^{i-1} W_{\widehat j}  = W_{[i,r]} \ .
\end{equation}
It follows that the weight $\gamma = u \omega_i$ satisfies 
both conditions in (2), so we must have $p_{u \omega_i} (x) = 0$.
But this contradicts the last statement in Lemma~\ref{lem:who-vanishes},
and we are done. 
\endproof

Notice that condition (2) in Proposition~\ref{prop:1} 
depends on the choice of ordering of fundamental weights.
We will introduce a special class of \emph{economical}
orderings that lead to the minimal possible number of 
equations in (2).

For any $i$, let $R(i)$ denote the set of positive roots 
whose expansion into the sum of simple roots contains 
the simple root~$\alpha_i\,$.

\begin{proposition}
\label{prop:2}
The correspondence $\alpha \mapsto s_\alpha \omega_i$
is an embedding of $R(i)$ into $W \omega_i - \{\omega_i\}$. 
\end{proposition}

\proof
Let $\alpha$ be a positive root. 
We have
\begin{equation}
\label{eq:reflections} 
\omega_i - s_\alpha \omega_i = (\omega_i, \alpha^\vee) \alpha \ ,
\end{equation}
where $(,)$ is a $W$-invariant scalar product of weights, and 
$\alpha^\vee = 2 \alpha / (\alpha, \alpha)$ is the dual root. 
By definition of fundamental weights, $(\omega_i, \alpha^\vee)$ is 
the coefficient of $\alpha_i^\vee$ in the expansion of $\alpha^\vee$ 
into the sum of dual simple roots.
Clearly this coefficient is nonzero precisely when $\alpha \in R(i)$.
Since no two positive roots are proportional to each other,
the vectors $(\omega_i, \alpha^\vee) \alpha$ for $\alpha \in R(i)$ 
are distinct nonzero vectors, proving the proposition.
\endproof

\begin{definition}
\label{def:economical-weight}
{\rm
We say that an index $i\in[1,r]$ (or the corresponding fundamental
weight~$\omega_i$) is \emph{economical} for $W$ if the correspondence
in Proposition~\ref{prop:2} is a \emph{bijection} between $R(i)$ and  
$W \omega_i - \{\omega_i\}$.
This is equivalent to 
\begin{equation}
\label{eq:econ}
1 + |R(i)| = |W \omega_i| = 
|W|/ |W_{\widehat i}|  \ .
\end{equation}
}\end{definition}

Here is a classification of all economical fundamental weights in 
irreducible Weyl groups.

\begin{proposition}
\label{prop:econ-classification}
Let $W$ be an irreducible Weyl group of rank $r$
with the set of simple reflections ordered as in \cite{bourbaki}.
An index $i$ is economical for $W$ precisely in the following three cases:

\noindent {\rm (1)} $r \leq 2$, and $i$ is arbitrary.

\noindent {\rm (2)} $W$ is of type $A_r$ for $r > 2$, 
and $i = 1$ or $i = r$. 

\noindent {\rm (3)} $W$ is of type $B_r$ or $C_r$ for $r > 2$, 
and $i = 1$. 

\end{proposition}  

\proof
First let us show that an index $i$ is indeed economical
in each of the cases (1) -- (3). 
The statement is trivial for type $A_1$.
If $r = 2$ then $W$ is of type $A_2$, $B_2$ or $G_2$,
i.e., is a dihedral group of cardinality $2d$ where $d = 3, 4$ or $6$,
respectively. 
We have $|W|/ |W_{\widehat i}| = d$ for any $i$ since $W_{\widehat i}$ is the 
two-element     group. 
On the other hand, $d$ is the number of positive roots in each case
which implies that $|R(i)| = d-1$ (the only positive root not in 
$R(i)$ is the simple root different from $\alpha_i$). 
Thus our statement follows from (\ref{eq:econ}). 

In case (2), we have $W = \mathcal{S}_{r+1}$ and 
$W_{\widehat i} = \mathcal{S}_{r}$ for $i = 1$ or $i = r$. 
Therefore, $|W|/ |W_{\widehat i}| = (r+1)!/r! = r+1 \,$. 
On the other hand, $R(1)$ (resp. $R(r)$) consists of $r$ roots 
$\varepsilon_1 - \varepsilon_{j+1}\,$
(resp. $\varepsilon_j - \varepsilon_{r+1}\,$) for $j = 1, \dots, r$, 
in standard notation of \cite{bourbaki}. 
Thus both $i = 1$ and $i = r$ are economical. 

Similarly, in case (3), the index  $i = 1$ (in the usual numeration)
is economical because 
$|W|/|W_{[2,r]}| = (2^r r!) / (2^{r-1} (r-1)!) = 2r$,
while $R(1)$ (say, for type $B_r$) consists of $2r-1$ roots:  
$\varepsilon_1 \pm \varepsilon_j$ ($j = 2, \dots, r$) and $\varepsilon_1$. 

To show that cases (1) -- (3) exhaust all economical indices,
we use the following observation: if $i$ is economical for 
$W$ then, in particular, we have
$$\omega_i - \wnot \omega_i = (\omega_i, \alpha^\vee) \alpha$$
for some positive root $\alpha$ (cf. (\ref{eq:reflections})),
where $\wnot$ is the maximal element of $W$. 
Since $\wnot$ sends positive roots to negative ones, 
it follows that $- \wnot \omega_i$ is also a fundamental weight
(possibly equal to $\omega_i$), and so $\alpha$ must be a dominant weight. 
If $W$ is simply-laced, i.e., all roots are of the same length,
then it is known that $W$ acts on the set $R$ of roots transitively.
Therefore, there is a unique root which is a dominant weight:
the maximal root $\alpha_{\max}$. 
The tables in \cite{bourbaki} show that if $W$ is 
simply-laced but not of type $A_r$ then $\alpha_{\max}$ is proportional 
to some fundamental weight $\omega_i$, so only this fundamental weight
has a chance to be economical.
But then we have  
$$|W \omega_i| = |W \alpha_{\max}| = |R| = 2 |R_+| > |R(i)| + 1 \ ,$$
so, for a simply laced $W$ not of type $A_r$, there are no economical
indices.

If $W$ is not simply-laced then there are precisely two 
roots which are dominant weights: the maximal long root and the maximal 
short root. 
Leaving aside cases (1) and (3) that we already considered,
this leaves only three more possibilites for an economical 
index: $i = 2$ for $W$ of type $B_r$
with $r > 2$; and $i = 1$ or $i = 4$ for $W$ of type $F_4$.
Since the root system of type $F_4$ is self-dual, we have
$|W \omega_1| = |W \omega_4| = |R_+|$,
while $|R(i)| \leq |R_+| - 3$ for any $i$ (since $R(i)$ does not contain
three simple roots different from $\alpha_i$). 
As for $W$ of type $B_r$ and $i = 2$, the set $W \omega_2$
consists of $2 r (r-1)$ weights of the form 
$\pm (\varepsilon_i \pm \varepsilon_j)$, $1 \leq i < j \leq r$,
and we have
$$|R(2)| + 1 \leq |R_+| - r + 2 = r(r-1) + 2  < 2 r (r-1) = |W \omega_2| 
\ .$$
We see that, in each of the three cases, $|R(i)| + 1 < |W \omega_i|$,
i.e., $i$ is not economical, and we are done. 
\endproof

\begin{proposition}
\label{prop:linear-ordering}
If a fundamental weight $\omega_i$ is economical for~$W$
then the Bruhat order on~$W\omega_i$ 
is \emph{linear}. 
\end{proposition}

\proof
Let $\gamma= w \omega_i$ and $\delta$ be two distinct 
Pl\"ucker weights of level~$i$.
Then $w^{-1} \delta \neq \omega_i\,$, 
which by Definition~\ref{def:economical-weight} implies that 
$w^{-1} \delta = t \omega_i$ for some reflection $t$. 
Since $wt$ and $w$ are comparable in the Bruhat order,
the same is true for $\delta = w t \omega_i$ and $\gamma= w \omega_i\,$.
\endproof

According to V.~Serganova (private communication), the converse of 
Proposition~\ref{prop:linear-ordering} is also true:
the Bruhat order on $W \omega_i$ is linear precisely in 
one of the cases (1)--(3) in Proposition~\ref{prop:econ-classification}.

\begin{definition}
\label{def:economical-ordering}
{\rm
A linear ordering of fundamental weights is called 
\emph{economical} if, for each~$i$, the index $i$ is economical for 
the group $W_{[i,r]}$. 
}
\end{definition}

This definition can be restated as follows. 
For a positive root $\alpha$, let $\mu (\alpha)$ 
denote the smallest index $i$ such that $\alpha \in R(i)$.
(In other words, the expansion of $\alpha$ does not contain the simple roots 
$\alpha_1,\dots,\alpha_{i-1}$ but does contain~$\alpha_i\,$.)
The ordering of fundamental weights is economical if and only if,
for every~$i\in[1,r]$, 
the map $\alpha\mapsto s_\alpha\omega_i$ is a bijection between 

\noindent
(i)\ \ the set of positive roots $\alpha$ with $\mu(\alpha)=i$ and

\noindent
(ii) the set $W_{[i,r]}\omega_i-\{\omega_i\}$.

Repeatedly using Proposition~\ref{prop:econ-classification},
we obtain the following corollary.

\begin{corollary}
\label{cor:standard econ}
An irreducible Weyl group possesses an economical ordering
of fundamental weights if and only if it is of one of the types 
$A_r, B_r, C_r$, or $G_2$.
In each of these cases, the standard ordering of fundamental weights 
given in~\cite{bourbaki} is economical. 
\end{corollary}

For an economical ordering, Proposition~\ref{prop:1} 
can be refined as follows. 

\begin{theorem}
\label{th:3}
Suppose the fundamental weights are ordered in an economical way.
Then an element $x \in G/B$ belongs to a Schubert cell 
$\cell_w$ if and only if:
\begin{align}
\label{eq:econ-neq}
&
\begin{array}{l}
\text{$p_{w \omega_i} (x) \neq 0$ for all $i$
such that there exists a positive root $\alpha$}\\
\text{with $\mu (\alpha) = i$
and $w \alpha$ negative;}
\end{array}
\\
\label{eq:econ-eq}
&\begin{array}{l}
\text{$p_{w s_\alpha \omega_{\mu (\alpha)}} (x) = 0$ 
for all positive roots $\alpha$ 
such that $w \alpha$ is also positive.}
\end{array}
\end{align}  
\end{theorem}

\proof
Recall that, for $\alpha>0$, the root 
$w \alpha$ is positive if and only if $w s_\alpha > w$. 
In view of this, Proposition~\ref{prop:1} shows that
conditions (\ref{eq:econ-neq})--(\ref{eq:econ-eq})
are indeed necessary. 

Assume that (\ref{eq:econ-neq})--(\ref{eq:econ-eq}) hold. 
To prove that $x\in\cell_w$, it suffices to show that 
$p_{w \omega_i} (x) \neq 0$ for \emph{all} $i\in [1,r]$. 
Suppose otherwise, and let $i$ be the minimal index such that 
$p_{w \omega_i} (x)=0$. 
By  (\ref{eq:econ-neq}), we have $w\alpha>0$ 
(thus $ws_\alpha>w$) for all positive roots $\alpha$
with $\mu(\alpha)=i$. 
In view of the definition of economical ordering, the weight
$w\omega_i$ is 
the minimal element of $wW_{[i,r]}\omega_i\,$. 
Now (\ref{eq:econ-eq}) implies that $p_\gamma(x)=0$ for all
$\gamma\in wW_{[i,r]}\omega_i-\{w\omega_i\}$. 

Suppose $x\in\cell_u\,$. 
The same argument as in the proof of Proposition~\ref{prop:1} shows
that $u\in wW_{[i,r]}$. 
Since $p_{u\omega_i}(x)\neq 0$, the weight $u\omega_i$ must coincide
with $w\omega_i\,$, which contradicts the assumption 
$p_{w \omega_i} (x)=0$. 
\endproof

The number of equations in (\ref{eq:econ-eq}) is equal to the number
of positive roots $\alpha$ such that $w \alpha$ is also positive;
this is precisely the codimension $\dim (G/B) - \l (w)$ of 
$\cell_w$ in the flag variety. 
Furthermore, the number of inequalities in (\ref{eq:econ-neq}) 
is at most $\min (r, \l (w))$. 
Applying Corollary~\ref{cor:standard econ}, we obtain the following 
solution of Problem~\ref{problem:short} for types 
$A$, $B$, $C$, and~$G_2$.

\begin{corollary}
\label{cor:ABCG}
For each of the types $A_r$, $B_r$, $C_r$, and $G_2$, conditions
{\rm (\ref{eq:econ-neq})--(\ref{eq:econ-eq})}
(with the standard ordering of fundamental weights)
describe an arbitrary Schubert cell $\cell_w$ using
$\dim (G/B) - \l (w)$ equations and at most $\min (r, \l (w))$~inequalities. 
\end{corollary}

As a special case, we obtain the following enhancement 
of \cite[Proposition~4.1]{FZ}. 

\begin{corollary}
\label{cor:Bruhat-SL}
For the type~$A_{n-1}$, 
an element $x \in G/B$ belongs to the Schubert cell $\cell_w$ if and only if
it satisfies the following conditions:
\begin{align}
\label{eq:rtl-minima}
&\text{$p_{w([1,i])}(x) \neq 0$ for all $i$
such that there exists $j>i$ with $w(j)<w(i)$;}\\
\label{eq:equations-A}
&\text{$p_{w([1,i-1] \cup \{j\})}(x) = 0$ whenever 
$1 \leq i < j \leq n$ and $w(i) < w(j)$.}
\end{align}
\nopagebreak
Thus $\cell_w$ can be described by at most $\binom{n}{2}$ equations
and inequalities of the form $p_I=0$ or $p_I\neq 0$. 
\end{corollary}

We conclude this section by addressing Problem~\ref{problem:short} for
type $D_r\,$. 
We  note that for $r\geq 4$, there are no economical indices. 
The index $i = 1$ (in the standard numeration) is ``one root short''
of being economical:  
$|W|/|W_{[2,r]}|=2r$ while $R(1)$ consists of $2r-2$ roots 
$\varepsilon_1 \pm \varepsilon_j$ ($j = 2, \dots, r$). 
As a consequence, we have to add extra equations to those in
(\ref{eq:econ-eq}) in order to describe~$\cell_w\,$. 
To minimize the number of these equations,
we use the following ordering of fundamental weights,
which is somewhat different from the one in~\cite{bourbaki}:

\begin{center}
\setlength{\unitlength}{1.5pt} 

\begin{picture}(105,50)(0,-5)

\put(0,20){\circle*{2.5}}
\put(20,20){\circle*{2.5}}
\put(60,20){\circle*{2.5}}
\put(80,20){\circle*{2.5}}
\put(100,0){\circle*{2.5}}
\put(100,40){\circle*{2.5}}

\put(0,20){\line(1,0){30}}
\put(80,20){\line(-1,0){30}}
\put(80,20){\line(1,1){20}}
\put(80,20){\line(1,-1){20}}

\put(40,19){\makebox(0,0){$\cdots$}}

\put(0,25){\makebox(0,0){$1$}}
\put(20,25){\makebox(0,0){$2$}}
\put(60,25){\makebox(0,0){$r-3$}}
\put(92,20){\makebox(0,0){$r-1$}}
\put(111,40){\makebox(0,0){$r-2$}}
\put(105,0){\makebox(0,0){$r$}}

\end{picture}
\end{center}

Theorem~\ref{th:3} and Corollary~\ref{cor:ABCG}
then have the following analogues (with similar proofs).

\begin{proposition}
\label{prop:D}
Let $G$ be of type $D_r$, $r\geq 4$,
and let the fundamental weights be ordered as above.
Then an element $x \in G/B$ belongs to a Schubert cell
$\cell_w$ if and only it satisfies conditions
{\rm (\ref{eq:econ-neq})--(\ref{eq:econ-eq})}, along with
the condition
\begin{equation}
\text{$p_\gamma (x) = 0$ 
whenever
$\gamma=w(\varepsilon_1\!+\cdots+\varepsilon_{i-1}\!-\varepsilon_i)
>w(\varepsilon_1\!+\cdots+\varepsilon_i)$, $i\leq r\!-\!3$.}
\end{equation}
Thus $\cell_w$ can be described using at most 
$\dim (G/B) - \l (w)+r-3$ equations and at most 
$\min (r, \l(w))$~inequalities.  
\end{proposition}

\section{Cell recognition algorithms}
\label{sec:recognition}

Our approach to the cell recognition problem 
(Problem~\ref{problem:recognition}) will be based on 
Proposition~\ref{prop:1} and Theorem~\ref{th:3}.

Suppose that the binary string $(b_\gamma)$ is the 
vanishing pattern of all Pl\"ucker coordinates at some 
point $x \in G/B$:
\begin{equation}
\label{eq:b-gamma}
b_\gamma = b_\gamma(x) = 
\begin{cases}
0 & \text{if $p_\gamma(x)=0$;}\\
1 & \text{if $p_\gamma(x)\neq0$.}
\end{cases}
\end{equation}

The following lemma is a reformulation of Lemma~\ref{lem:who-vanishes}.

\begin{lemma}
\label{lem:unique-max}
For any $x \in G/B$ and any $i\in[1,r]$, the set of all Pl\"ucker
weights $\gamma$ of level $i$ such that $b_\gamma (x)=1$ has a
unique maximal element with respect to the Bruhat order on $W
\omega_i$. 
Furthermore, if $x$ belongs to the Schubert cell 
$\cell_w = (BwB)/B$, 
then this maximal element is equal to~$w \omega_i$.
\end{lemma}

In view of Lemma~\ref{lem:unique-max}, any vector $b_\gamma(x)$ is
``acceptable'' according to the following definition. 

\begin{definition}
\label{def:acceptable}
{\rm
A binary vector $(b_\gamma)$, where $\gamma$ runs over all  Pl\"ucker
weights, is called \emph{acceptable} if
\begin{align}
\label{eq:accept1}
&
\begin{array}{l}
\text{for any $i\in[1,r]$, the set $\{\gamma\in W\omega_i \,:\,
  b_\gamma=1\}$ is nonempty, and has  
}\\
\text{a unique maximal element $\gamma_i$ with respect to the
  Bruhat order;}
\end{array}
\\
\label{eq:accept2}
&\begin{array}{l}
\text{there exists $w\in W$ such that $\gamma_i=w\omega_i$ for any~$i$.}
\end{array}
\end{align}  
}\end{definition}

It is immediate from Lemma~\ref{lem:deodhar}
that the element $w$ in (\ref{eq:accept2}) is unique.

We will now study the following purely combinatorial problem that
includes Problem~\ref{problem:recognition}  as a special case. 

\begin{problem}
\label{problem:acceptable}
{\rm
For a given acceptable vector $(b_\gamma)$, 
compute the element~$w$ in (\ref{eq:accept2}) 
by testing the minimal number of bits~$b_\gamma$. 
}\end{problem}
 
For $\gamma \in W \omega_i\,$, let us denote 
$W(\gamma) = \{u \in W: u \omega_i = \gamma\}$. 
Thus $W(\gamma)$ is a left coset in $W$ with respect to the 
stabilizer of $\omega_i$ (i.e., with respect to $W_{\widehat{i}}$). 
Our approach to Problem~\ref{problem:acceptable} will be based on the
following lemma, which follows
from~(\ref{eq:intersection-of-parabolics}).  

\begin{lemma}
\label{lem:accept-properties}
Let $(b_\gamma)$ be an acceptable binary vector. 
In the notation of Definition~\ref{def:acceptable}, for every $i$, we
have: 
\[
W(\gamma_1) \cap \cdots \cap W(\gamma_{i-1}) 
= w W_{[i,r]}\ ;
\]
also, 
$\gamma_i$ is the maximal element of 
$ w W_{[i,r]}\omega_i$ such that $b_{\gamma_i} =1$.
\end{lemma}

The following algorithm for Problem~\ref{problem:acceptable} 
is based on Lemma~\ref{lem:accept-properties}; 
it successively computes the weights  
$\gamma_1, \gamma_2, \dots$, 
and in the end obtains $w$ as the sole element in the intersection
$W(\gamma_1) \cap \cdots \cap W(\gamma_{r})$.

\begin{algorithm}
\label{alg:recognition-general}
{\rm
 ~ \\
\begin{tabular}{ll}
\emph{Input:} 
& acceptable binary vector $(b_\gamma)$.\\ 
\emph{Output:} 
& the element $w\in W$ given by (\ref{eq:accept2}). 
\end{tabular}
\nopagebreak

\begin{tabular}{l}
$U:=W$;\\
\textbf{for $i$ from $1$ to $r$ do}\\
\quad fix a linear order $U\omega_i\!=\!\{\eta_1<\cdots<\eta_m\}$ 
compatible with the Bruhat
         order;\\ 
\quad $j:=m$; \\
\quad \textbf{while $b_{\eta_j}=0$ do $j:=j-1${\rm ;} od};\\
\quad \textbf{comment: } 
  $\eta_j=\gamma_i=\max\{\gamma\in U\omega_i\,:\, b_\gamma=1\}$\\
\quad $U:=U\cap W(\eta_j)$;\\
\textbf{od};\\
\textbf{return}(U);
\end{tabular}
}
\end{algorithm}

In particular, this algorithm can be used to solve
Problem~\ref{problem:recognition}:
if the input vector $(b_\gamma)$ is the vanishing pattern
(\ref{eq:b-gamma}) for a point $x\in G/B$, 
then the algorithm returns the element~$w\in W$ such
that~$x\in\cell_w$. 

The algorithm depends on the choice of the ordering of fundamental
weights. As in Section~\ref{sec:short}, the best results are achieved
for economical orderings. 
In this case, Proposition~\ref{prop:linear-ordering} implies that the
set of weights 
$
U\omega_i=w W_{[i,r]} \omega_i
$
appearing in Algorithm~\ref{alg:recognition-general} is linearly
ordered by the Bruhat order, making the third line of the algorithm
redundant. 

In particular, in the case of type~$A_{n-1}$, 
the standard ordering of the fundamental weights, 
and an acceptable vector defined by~(\ref{eq:b-gamma}), 
Algorithm~\ref{alg:recognition-general}
takes the following form. 
(As before, we identify the Pl\"ucker weights 
with subsets in $[1,n]$.)

\begin{algorithm}
\label{alg:recognition-gln}
{\rm
 ~ \\
\begin{tabular}{ll}
\emph{Input:} 
& vanishing pattern of Pl\"ucker coordinates of a complete
flag~$x$ in~$\CC^n$.\\ 
\emph{Output:} 
& permutation $w\in \mathcal{S}_n$ such that $x\in\cell_w\,$. 
\end{tabular}
\nopagebreak

\begin{tabular}{l}
$I:=\emptyset$;\\
\textbf{for $i$ from $1$ to $n$ do}\\
\quad $k:=n$;\\
\quad \textbf{while $k>\min([1,n]\!-\!I)$ 
               and $(k\in I$ or $p_{I\cup\{k\}}(x)=0)$ do $k:=k\!-\!1$; od};\\
\quad $w(i):=k$;\\
\quad $I:=I\cup\{k\}$;\\
\quad \textbf{comment:} $I=w([1,i])$\\
\textbf{od};
\end{tabular}
}
\end{algorithm}

To convince oneself that Algorithm~\ref{alg:recognition-gln} is a
specialization of Algorithm~\ref{alg:recognition-general},
it suffices to observe the following:
the weights in $wW_{[i,r]}\omega_i$ correspond to 
the $i$-subsets of the form $w([1,i-1])\cup \{k\}$, 
and the Bruhat order on $wW_{[i,r]}\omega_i$ 
corresponds to the usual ordering of the values~$k$. 

In the special case of type~$A_2\,$, we recover the algorithm
presented in Figure~\ref{fig:3steps}. 

Algorithm~\ref{alg:recognition-gln} agrees completely 
with the description of Schubert cells
given in Corollary~\ref{cor:Bruhat-SL}: 
to arrive at any $w$, we need to check 
exactly the same Pl\"ucker coordinates that appear in 
(\ref{eq:rtl-minima})--(\ref{eq:equations-A}). 
We thus obtain the following result.

\begin{proposition}
\label{prop:recognition-gln}
For a complete flag $x$ in $\CC^n$, 
Algorithm~\ref{alg:recognition-gln} recognizes the Schubert cell $x$
is in by testing at most $\binom{n}{2}$ bits of the vanishing pattern
of its Pl\"ucker coordinates. 
\end{proposition}

We omit the type $B$ (or~$C$) analogues of
Algorithm~\ref{alg:recognition-gln} 
and Proposition~\ref{prop:recognition-gln},
which can be obtained in a straightforward way. 

\section{
On the number of equations
defining a Schubert variety}
\label{sec:lower-bounds-for-the-number-of-equations}

Problem~\ref{problem:short} is closely related to the classical
problem of describing Schubert varieties $X_w$ 
as algebraic subsets of~$G/B$.

\begin{problem}
\label{problem:variety description}
\emph{(Short descriptions of Schubert varieties)}
Define an arbitrary Schubert variety $X_w$ 
(as a subset of $G/B$) by as small as possible 
number of equations of the form $p_\gamma =0$. 
\end{problem}

The aim of this section is to demonstrate that,
for a certain Schubert variety~$X_w$ of type~$A_{n-1}\,$,
one needs exponentially many (as a function of~$n$)
such equations to define $X_w$ (set-theoretically).

Throughout this section, $G=SL_n$ and $W=\mathcal{S}_n\,$.
Any Schubert cell $\cell_w$ has the  special representative 
$\pi_w$: it is a complete flag in $\CC^n$ formed by 
the coordinate subspaces $\CC^{w([1,i])}$ for $i = 1, \dots, n$.  
The following obvious observation will be useful in
obtaining lower bounds. 

\begin{lemma}
\label{lem:perm-matrices} 
For $w\in\mathcal{S}_n\,$, 
a Pl\"ucker coordinate $p_I$ does not vanish at $\pi_w$ 
if and only if $I=w([1,|I|])$. 
\end{lemma}

\begin{proposition}
\label{prop:lower-bound}
Suppose that $n = 4k$ is divisible by $4$. 
Let $w\in\mathcal{S}_n$ be the maximal element 
of the parabolic subgroup 
$W_{\widehat {2k}} = \mathcal{S}_{2k} \times \mathcal{S}_{2k} \subset
\mathcal{S}_n$
(thus $w$ puts the elements in each of the blocks $[1,2k]$ and
$[2k+1, 4k]$ in the reverse order).
Suppose the set $\mathcal{I}$ is such that
\[
X_w=\{x\in G/B\,:\,p_I(x)=0 \text{ for } I \in \mathcal{I}\}.
\]
Then 
\begin{equation}
\label{eq:lower-bound}
|\mathcal{I}| \geq  \binom{2k}{k} \ .
\end{equation}
\end{proposition}

Note that the right-hand side of (\ref{eq:lower-bound})
grows as $2^{n/2}/\sqrt{n}$,
while the codimension of this particular Schubert variety $X_w$
equals $(n/2)^2$. 

\proof
Our lower bound for $|\mathcal{I}|$ is based on the following idea.
Suppose a permutation $u \in \mathcal{S}_n$ is such that 
$u \not\leq w$.
Then the flag $\pi_u$ does not belong to the Schubert variety 
$X_w$, so there must exist $I \in \mathcal{I}$ such that 
$p_I (\pi_u) \neq 0$.
By Lemma~\ref{lem:perm-matrices}, this means that $I=u([1,|I|])$.
In view of Lemma~\ref{lem:who-vanishes}, the inclusion 
$I \in \mathcal{I}$ also implies that 
$I \not\leq w([1,|I|])$.  
We conclude that, in order to prove (\ref{eq:lower-bound}),
it suffices to construct a subset $U \subset \mathcal{S}_n$
satisfying the following three properties:

\noindent
(1) $u \not\leq w$ for any $u \in U$;

\noindent
(2) $|U| = \binom{2k}{k}^2$;

\noindent
(3) for every subset $I \subset [1,n]$ such that 
$I \not\leq w([1,|I|])$, there are at most $\binom{2k}{k}$
permutations $u \in U$ such that $I = u([1,|I|])$. 

Define $U$ to be the set of all permutations $u$ that send 
$[1,k] \cup [2k+1,3k]$ onto $[1,2k]$, and increase 
on each of the blocks $[1,k]$, $[k+1,2k]$, $[2k+1,3k]$, and $[3k+1,4k]$. 
Each $u \in U$ is uniquely determined by two $k$-subsets $A = u([1,k])
\subset [1,2k]$ and $B = u([k+1,2k]) \subset [2k+1,4k]$;
we write $u=u_{A,B}$. 
Now (2) is obvious.  
Since $u_{A,B} ([1,2k]) = A \cup B > [1,2k] = w([1,2k])$, 
we have $u \not\leq w$ for any $u \in U$, so $U$ satisfies (1). 

It remains to prove (3).
Let $I \subset [1,n]$ be such that $I \not\leq w([1,|I|])$.
We need to show  that there are at most $\binom{2k}{k}$
permutations $u_{A,B} \in U$ such that $I = u_{A,B} ([1,|I|])$. 
First of all, we have $u_{A,B} ([1,i]) \leq w([1,i])$
for $i \leq k$ or $i \geq 3k$. 
Therefore, we may assume that $k < |I| < 3k$. 
Let us consider two cases.

\noindent {\bf Case 1.} $|I| = k+l$ for some $l \in [1,k]$. 
The equality $I = u_{A,B} ([1,|I|])$ means that $I$ is the union 
of $A$ and the set of $l$ smallest elements of $B$. 
Thus $A = [1,2k] \cap I$ is uniquely determined by $I$, 
while the number of choices for $B$ is $\binom{4k - \max I}{k-l}$,
which is less than $\binom{2k}{k}$. 

\noindent {\bf Case 2.} $|I| = 2k+l$ for some $l \in [1,k-1]$. 
Now the equality $I = u_{A,B} ([1,|I|])$ means that $I$ is the union 
of $A$, $B$, and the set of $l$ smallest elements of $[1,2k] - A$. 
Thus $B = [2k+1,4k] \cap I$ is uniquely determined by $I$, 
while the number of choices for $A$ is
$\binom{k+l}{k} < \binom{2k}{k}$. 

This concludes the proof of (\ref{eq:lower-bound}). 
\endproof
                 
\begin{corollary}
There exist elements $u < v$ in $W=\mathcal{S}_{4k}$ such that 
$X_u$ has codimension~1
in~$X_v\,$, while defining $X_u$ inside~$X_v$ requires at least 
$\displaystyle\frac{1}{4k^2}\binom{2k}{k}$ equations of the form~$p_I=0$. 
\end{corollary}

\proof
Consider a saturated chain $w=v_0<v_1<\cdots<v_N=\wnot$ in the Bruhat
order, where $w$ is the same as in Proposition~\ref{prop:lower-bound}.
(thus $N= 4k^2$).
If $M(u,v)$ denotes the minimal number of equations of the
form~$p_I=0$ defining $X_u$ inside $X_v\,$, then obviously
$M(w,\wnot)\leq \sum M(v_i,v_{i+1})\leq N\cdot\max_i(M(v_i,v_{i+1}))$. 
Combining this with the lower bound on $M(w,\wnot)$ obtained in 
Proposition~\ref{prop:lower-bound} completes the proof. 
\endproof

\section{On cell recognition without feedback}
\label{sec:static}

In this section, we examine the following problem. 

\begin{problem}
\label{problem:static}
\emph{(Cell recognition without feedback)}
Find a subset of Pl\"ucker coordinates of smallest possible
cardinality whose vanishing pattern at any point $x\in G/B$ uniquely
determines the Schubert cell of~$x$.
\end{problem} 

Notice that, unlike in Problem~\ref{problem:short}, the Schubert cell
is not fixed in advance; and in contrast to
Problem~\ref{problem:recognition}, we have to present the entire list
of Pl\"ucker coordinates right away (i.e., there is no feedback). 

\begin{example}
\label{example:gl3-static}
{\rm
Consider the special case of $G=SL_3\,$. 
Analyzing Table~\ref{table:gl3} in Section~\ref{sec:gl3},
we discover that the list in question must contain the Pl\"ucker
coordinates~$p_3$ (to distinguish between vanishing patterns of generic
elements of Schubert cells labelled by $s_1s_2$ and~$\wnot$),
$p_2$ (same reason, for $e$ and~$s_1$), $p_{13}$ (for $e$ and~$s_2$),
and $p_{23}$ (for $s_2s_1$ and~$\wnot$). 
The vanishing pattern of these 4 Pl\"ucker coordinates does indeed
determine the cell a point is in (see last column of
Table~\ref{table:gl3}).  
Hence this 4-element collection of Pl\"ucker coordinates provides the
unique solution 
to Problem~\ref{problem:static} for the type~$A_2$. 
}\end{example}

The following result shows that for the type~$A$, 
the subset asked for in Problem~\ref{problem:static} must contain an
overwhelming proportion of all Pl\"ucker coordinates. 

\begin{proposition}
\label{prop:static}
For the type~$A_{n-1}$,
any subset satisfying the requirements in Problem~\ref{problem:static}
contains at least the $\frac{n-1}{n+1}$ proportion of all Pl\"ucker
coordinates. 
\end{proposition} 

Note that there are $2^n-2$ Pl\"ucker coordinates altogether
in this case.  

\proof
We will actually show more: that this many Pl\"ucker coordinates are
needed to distinguish between the vanishing patterns of any two
different elements of the form $\pi_w$, for $w\in W=\mathcal{S}_n$  
(we use the notation introduced at the beginning of
Section~\ref{sec:lower-bounds-for-the-number-of-equations}).
Let $\mathcal{I}$ be a collection of subsets $I\subset[1,n]$ such that
the vanishing patterns of the Pl\"ucker coordinates $p_I(\pi_w)$, for
$I\in\mathcal{I}$, are distinct for all elements~$w\in W$. 
In view of Lemma~\ref{lem:perm-matrices}, this means that for any
distinct $u,v\in W$, there exists an index~$i\in [1,n]$ such that
the subsets $u([1,i])$ and $v([1,i])$ are distinct,
and at least one of them belongs to~$\mathcal{I}$. 

Let $I$ be a nonempty proper subset of~$[1,n]$ of cardinality~$i$. 
Choose $u\in W$ so that $u([1,i])=I$, and let $v=us_i\,$. 
Then $u([1,j])=v([1,j])$ unless~$j=i$, implying that 
$\mathcal{I}$ must contain either $u([1,i])=I$ or 
$v([1,i])=I\setminus \{u(i)\} \cup \{u(i\!+\!1)\}$ (or both). 
We conclude that for any two subsets $I,J\subset[1,n]$ of the same
cardinality which are Hamming distance~2 from each other
(i.e., one is obtained from another by exchanging a single element), 
the collection $\mathcal{I}$ has to contain either $I$ or~$J$. 

Let $\overline{\mathcal{I}}_i$ denote the collection of all
$i$-subsets of $[1,n]$ not in~$\mathcal{I}$. 
Then $\overline{\mathcal{I}}_i$ does not contain two subsets at
Hamming distance~2 from each other. 
Such collections of subsets are called \emph{binary codes of constant
weight detecting single errors},
and they were an object of extensive study in coding theory. 
In particular, various upper bounds on the cardinality of such a code
have been obtained; see, for example,
\cite[Chapter~17]{macwilliams-sloane}. 
(We thank Richard Stanley for providing this reference.) 
For our purposes, it will suffice to have a very simple upper bound 
\begin{equation}
\label{eq:code-bound}
|\overline{\mathcal{I}}_i|\leq \frac{1}{i}\binom{n}{i-1} 
= \frac{1}{n+1}\binom{n+1}{i}\ .
\end{equation}
Although this bound is immediate from a sharper 
\cite[Ch.~17, Corollary~5]{macwilliams-sloane}, we will give a proof  
for the sake of completeness. 
 
To prove (\ref{eq:code-bound}), note that all $(i-1)$-subsets
contained in 
various $i$-subsets in $\overline{\mathcal{I}}_i$ must be distinct. 
Each $I\in \overline{\mathcal{I}}_i$ contains $i$ such subsets,
implying that $i\cdot |\overline{\mathcal{I}}_i|\leq \binom{n}{i-1}$,
as desired. 

The proof of Proposition~\ref{prop:static} can now be completed as
follows: 
\[
\begin{array}{rcl}
|\mathcal{I}|
&=& 2^n-2-\sum_{i=1}^{n-1}|\overline{\mathcal{I}}_i| \\[.1in]
&\geq& 2^n-2-\frac{1}{n+1}\sum_{i=1}^{n-1} \binom{n+1}{i}  \\[.1in]
&=& 2^n-2-\frac{1}{n+1}(2^{n+1}-n-3) \\[.1in]
&=& \frac{n-1}{n+1}(2^n-1)\ .\endproofmath \\[.1in]
\end{array}
\]

\section{Generic vanishing patterns 
}
\label{sec:bigrassmannian}

In the course of the above proof of Proposition~\ref{prop:static}, 
we have actually shown the following: 
assuming there is no feedback, 
``almost all'' Pl\"ucker coordinates are needed 
to distinguish between special representatives 
$\pi_w$ of Schubert cells. 
We will now demonstrate that the situation changes dramatically if we
replace these ``most special" representatives by the 
``most generic" ones.

In what follows, $W$ is an arbitrary Weyl group.
We associate to any $w \in W$ 
the \emph{generic} vanishing pattern $(b^{\rm gen}_\gamma (w))$ 
defined by
\begin{equation}
\label{eq:b-generic}
b^{\rm gen}_\gamma (w) =
\begin{cases}
1 & \text{if $\gamma \leq w \omega_i\,$;}\\
0 & \text{if $\gamma \not\leq w \omega_i\,$,}
\end{cases}
\end{equation}
where $\gamma$ runs over all Pl\"ucker weights of any level~$i$.
By Lemma~\ref{lem:who-vanishes}, this is the vanishing pattern 
$(b_\gamma (x))$ (cf. (\ref{eq:b-gamma})) of Pl\"ucker coordinates
for a generic element $x \in \cell_w$. 
                         
\begin{problem}
\label{problem:generic}
\emph{(Recognizing generic points without feedback)} 
Find a minimal subset of Pl\"ucker coordinates 
whose vanishing pattern distinguishes between 
the generic patterns $(b^{\rm gen}_\gamma (w))$.
\end{problem} 

Our solution of this problem will be based on the techniques
developed by Lascoux and Sch\"utzenberger~\cite{Las},
and further enhanced by Geck and Kim~\cite{Geck-Kim}. 
Let us first recall the main definitions and results of these papers. 

Let $P$ be a finite poset with unique minimal and maximal elements.
We say that $a\in P$ is the \emph{supremum} of a subset $Q\subset P$ if  
$a\geq q$ for any $q\in Q$, and moreover $a < b$ 
for any other element $b \in P$ with this property. 

\begin{definition}{\rm
The \emph{base} $B=B(P)$ of $P$ is the subset of $P$ consisting of all 
elements $a\in P$ which cannot be obtained as the supremum of a
subset of~$P$ not containing~$a$. 
}
\end{definition}

\begin{proposition}
\label{prop:Las} \cite{Las}
The map $a\mapsto \{b\in B\,:\,b\leq a\}$
is an embedding of $P$ (as an induced subposet) into the boolean
algebra of all subsets of~$B=B(P)$. 
Moreover, any other subset $B'\subset P$ with this property
contains~$B$. 
\end{proposition}

The following result appeared in \cite[Th\'eor\`eme~3.6]{Las};
another proof was given in \cite[Theorem~2.5]{Geck-Kim}. 

\begin{theorem}
\label{th:Las} 
\cite{Las}
For every element $u$ in the base of a finite Coxeter group~$W$, 
there are unique simple reflections $s_i$ and $s_j$ 
such that $us_i<u$ and $s_j u <u$. 
\end{theorem}

Let $\mathfrak{B}(W)$ denote the subset of Pl\"ucker weights which
correspond to the elements of the base $B(W)$, as follows: 
\[
\mathfrak{B}(W) = 
\{ u\omega_i \,:\, u\in B(W)\,,\  us_i<u \} \,.
\]

\begin{proposition}
\label{prop:bigrassmannian}
The correspondence $w \mapsto (b^{\rm gen}_\gamma (w))$,
where $\gamma$ runs over $\mathfrak{B}(W)$, 
is an embedding of $W$ (as an induced subposet) into 
the Boolean lattice of all binary vectors of the corresponding
length. 
Moreover, $\mathfrak{B}(W)$ is a minimal subset of Pl\"ucker weights 
that has this property.  
\end{proposition}

Thus the set of the Pl\"ucker coordinates
$p_\gamma\,$, with $\gamma\in\mathfrak{B}(W)$, 
provides a solution of Problem~\ref{problem:generic}.

\proof
Let $u\in B(W)$, and let $\gamma = u \omega_i\in\mathfrak{B}(W)$ 
be the corresponding weight. 
Since $u$ is the minimal representative of the coset 
$uW_{\widehat i}$, it follows that for any $w\in W$, the condition
``$\gamma \leq w \omega_i$'' 
is equivalent to ``$u \leq w$.''
Therefore, (\ref{eq:b-generic}) becomes
\begin{equation}
\label{eq:b-generic-bigrassmannian}
b^{\rm gen}_\gamma (w) =
\begin{cases}
1 & \text{if $u \leq w$;}\\
0 & \text{if $u \not\leq w$.}
\end{cases}
\end{equation}
Thus the set of non-vanishing Pl\"ucker coordinates $p_\gamma$,
$\gamma\in\mathfrak{B}(W)$, at a generic point 
in $\cell_w$ corresponds exactly to the set of elements 
in the base $B(W)$ 
that are less than or equal than~$w$ in the Bruhat order. 
The proposition then follows from Proposition~\ref{prop:Las}. 
\endproof  

The bases $B(W)$ were explicitly described and enumerated
in~\cite{Las} (for the types $A$ and~$B$) and~\cite{Geck-Kim} (for all
other types). 
As shown in~\cite{Geck-Kim,Las}, if $W$ is of one of the classical types
$A_r$, $B_r$, and $D_r$, then the cardinality of $B(W)$ is a cubic
polynomial in~$r$. 
In particular, for the type $A_{n-1}$ 
when $W = \mathcal{S}_n$, the base consists of the $\binom{n+1}{3}$
``bigrassmannian'' permutations: every 
triple of integers $0 \leq a < b < c \leq n$ gives rise to 
a such a permutation that acts identically 
on each of the blocks $[1,a]$ and $[c+1,n]$ while interchanging 
the blocks $[a+1,b]$ and $[b+1,c]$. 
The corresponding \emph{bigrassmannian Pl\"ucker coordinate} is 
$p_{[1,a]\cup [b+1,c]}$. 
Proposition~\ref{prop:bigrassmannian} tells that the vanishing pattern
of these $\binom{n+1}{3}$ Pl\"ucker coordinates
uniquely determines the Schubert cell of a given complete flag
$x$ in~$\CC^n$, provided we know that $x$ is generic within its
cell.  
In the special case $n=3$, the bigrassmannian Pl\"ucker coordinates
are exactly the four coordinates $p_2,p_3,p_{13},p_{23}$ 
involved in Example~\ref{example:gl3-static}
and in the descriptions of Section~\ref{sec:gl3}.

\section*{Acknowledgments}

We thank V.~Lakshmibai, Alain Lascoux, Peter Littelmann, Peter Magyar,
and Richard Stanley for helpful comments.

\nopagebreak

\end{document}